\input amstex
\documentstyle{amsppt}
\magnification=1200
\pageheight{43pc}
\pagewidth{29.5pc}
\hoffset=7mm
\voffset=1.2cm
\nopagenumbers
\NoRunningHeads
\NoBlackBoxes
\nologo
\topmatter
\title\nofrills ON DIRECTED ZERO-DIVISOR GRAPHS OF FINITE RINGS
\endtitle
\thanks
\endthanks
\author Tongsuo Wu
\endauthor
\email wutsc\@online.sh.cn
\endemail
\affil Department of Mathematics,
Shanghai Jiaotong University,
Shanghai 200030, P. R. China
\endaffil
\date{}\enddate
\abstract For an artinian ring $R$, the directed zero-divisor
graph $\Gamma(R)$ is connected if and only if there is no proper
one-sided identity element in $R$. Sinks and sources are
characterized and clarified for finite ring $R$, especially, it is
proved that for {\it any} ring $R$, if there exists a source $b$
in $\Gamma(R)$ with $b^2=0$, then $|R|=4$ and $R=\{0,a,b,c\}$,
where $a$ and $c$ are left identity elements and $ba=0=bc$. Such a
ring $R$ is also the only ring such that $\Gamma (R) $ has exactly
one source. This shows that $\Gamma(R)$ can not be a network for
any ring $R$.
\endabstract
\endtopmatter
\document
\baselineskip=14pt
\head 1. Introduction\endhead
\vskip3mm
  For any noncommutative ring $R$, let $Z(R)$ be
the set of (one-sided) zero-divisors of $R$. The directed
zero-divisor graph of $R$ is a directed graph $\Gamma(R)$ with
vertex set $Z(R)^*=Z(R)-\{0\}$, where for distinct vertices $x$
and $y$ of $Z(R)^*$ there is a directed edge from $x$ to $y$ if
and only if $xy=0$ ([6]). This is a generalization of zero-divisor
graph of commutative rings. The concept of a zero-divisor graph of
a commutative ring was introduced in [4], and it was mainly
concerned with colorings of rings there. In [2], the vertex set of
$\Gamma(R)$ is chosen to be $Z(R)^*$ and the authors study
interplay between the ring-theoretic properties of commutative
ring $R$ and the graph-theory property of $\Gamma(R)$. The
zero-divisor graph of a commutative ring is also studied by
several other authors, see [3] for a list of references. The
zero-divisor graph has been also introduced and studied for
semigroups in [5].\par In this paper, we study the directed
zero-divisor graph of noncommutative rings and we focus our
attention on finite rings (Most results on finite rings in this
paper actually holds for artinian rings.). In section 2, we prove
that an artinian ring $R$ has connected zero-divisor graph if and
only if one-sided identity of $R$ (if exists) is two-sided, if and
only if $\Gamma(R)$ contains no end vertex (i.e., sinks and
sources). For any distinct vertices $x,y$ of a finite ring with
proper one-sided identity, the directed distance from $x$ to $y$
is less than 7, if a directed path exists from $x$ to $y$. In
section 3 and 4, we study sinks and sources of finite and infinite
ring. We proved that for {\it any} ring $R$, if there exists a
source $b$ in $\Gamma(R)$ with $b^2=0$, then $|R|=4$ and
$R=\{0,a,b,c\}$, where $a$ and $c$ are left identity elements and
$ba=0=bc$. Such a ring $R$ is also the only ring such that $\Gamma
(R) $ has exactly one source. The dual result for sink vertex is
also true. This result is key to clarify sinks and sources in
$\Gamma(R)$. In Section 3, we show that for a finite ring $R$ with
at least five elements, sink and source can not coexist in
$\Gamma(R)$. For finite ring $R$ with at least five elements,
sinks (sources) are characterized by strongly right (left)
invertible elements. In section 4, we show that for any ring $R$
with at least five elements, there are only four possibilities for
the semigroups (when non-empty) $Sink(R)$, the set of sinks of
$\Gamma(R)$, and $Sour(R)$, the semigroup of all sources:\par (1)
$Sink(R)=\varnothing$ and $Sour(R)=\varnothing$;\par (2)
$Sink(R)=\varnothing$ and $|Sour(R)|\geq 2$;\par (3)
$Sour(R)=\varnothing$ and $|Sink(R)|\geq 2$;\par (4) $
|Sour(R)|=\infty =|Sink(R)|$.\par \noindent Therefore $\Gamma(R)$
can not be a network for any ring $R$.\par \vskip3mm \head 2.
Connectedness and Diameters \endhead \vskip3mm

\proclaim{Lemma 2.1}
 Suppose that every right identity element of a
finite ring R is a two-sided identity. We have \par
(1) Each left zero-divisor is also a right zero-divisor;\par
(2) If in addition $|R|\geq 5$, then for distinct
$a,b\in R^*$ with $ab=0$, there exists $c\in R^*$ such that
$ca=0$ and $c\neq a$.\par
\endproclaim
\vskip3mm
\demo{Proof}
(1) If $ab= 0$ and $xa\neq 0$ for all $x\in R$, then
$Ra=R$ and there is a right identity element $e\in R$.
Then $e$ is
also a left identity. Let $ga=e$,then $0=gab=eb=b$, a
contradiction.\par
(2) Suppose that $ann_l(a)=\{a,0\}$, then for distinct
$b,c,d\in R^*-\{a\}$, we have $daa=0,da\neq0$.
Thus $da=a.$ Then we have $b-c=b-d$ and therefore,
$c=d$. Then $|R|\leq 4$. This completes the proof.
\hfill QED
\enddemo
\vskip3mm
\proclaim{Lemma 2.2} Suppose $|R|\leq 4$ and every
proper one-sided identity element is a two-sided identity. Then
for distinct $a,b\in R^*$ with $ab=0$, we have $ba=0$.
\endproclaim
\vskip3mm
\demo{Proof}
If $R$ has three elements, then for $a\neq b$ with $ab=0$,
one has $(a+b)b=0$. Thus $ba=b(a+b)=0$.\par

Now suppose that $R=\{0,a,b,a+b\}$ has four elements. Without
loss, we assume $ab=0$. If $R$ has identity element $1$, assume
$a+b=1$. Then $a^2=a,b^2=b$. So if $ba=1$, then $a=aba=0$; if
$ba=a$, then $a=a^2=aba=0$; If $ba=b$, then $b=b^2=bab=0$. There
is contradiction in each case. Hence $ba=0$. In the remaining part
of the proof, assume that $R$ has no one-sided identity element.
Without loss, we can assume $2a=2b=0$ (The only other case is that
the additive group of $R$ is cyclic of order four.). Now we show
that $ab=0$ implies $ba=0$:\par \indent (1)$ba\neq a$: If $ba=a$,
then $a^2=0.$ In this case, we assert $b^2=0$. Actually, $b^2\neq
b$ since otherwise $b$ is a left identity element of $R$; $b^2\neq
a$ since otherwise $a=ba=b^2a=a^2=0$; $b^2\neq a+b$ since
otherwise, $b^2=(a+b)b=b^3=b(a+b)=a+(a+b)=b$. But $b^2=0$ implies
$a=ba=b^2a=0$, a contradiction.\par \indent (2)$ba\neq b$: If
$ba=b$, then $b^2=0.$ Then $a^2\neq a$ since otherwise $a(a+b)=a,$
$b(a+b)=b$, $(a+b)^2=a+b$, i.e., $a+b$ is a right identity element
of $R$; $a^2\neq 0$ since $b=ba^2$;$a^2\neq b$ since otherwise
$b=ba=ba^2=b^2=0$. Finally, $a^2=a+b$ and we have
$a^2=a+b=a^3=a^2+ba=a^2+b$, contradicting with the assumption
$b\neq 0$.\par \indent (3) $ba\neq a+b$: If $ba=a+b$, then
$a^2=0$, $0=ba^2=a^2+ba$, $0=ba=a+b$, another contradiction. This
completes the proof of $ba=0$, for rings without one-sided
identity element. \hfill QED
\enddemo
\vskip3mm Let $K_i$ be the complete directed graph with $i$
vertices. For any ring $R$, by [6, Theorem 3.2], there is no
isolated vertices in $\Gamma(R)$. Thus by Lemma 2.2, we have a
list of all possibilities of $\Gamma(R)$ for rings $R$ with
$|R|\leq 4$:\par (1) $|R|=2$: $K_i,i=0,1$\par (2) $|R|=3$:
$K_i,i=0,2$\par (3) $|R|=4$: $K_i$ ($i=0,1,2,3$);
$\circ\leftrightarrows\circ\leftrightarrows\circ$;
$\circ\rightarrow \circ\leftarrow$;
$\circ\leftarrow\circ\rightarrow\circ $. \par We remark that every
graph in the list can be realized as the zero-divisor graph of
some ring $R$ with $|R|\leq 4$. Actually this list is a special
case of results in [7, Section 4])\par
 \vskip3mm
\proclaim{Proposition 2.3} Suppose that $R$ is a finite ring with
the property that every one-sided identity element is a two-sided
identity in $R$. Then for any path $a\rightarrow b$ in
$\Gamma(R)$, there is a walk $c\rightarrow a\rightarrow
b\rightarrow d$, where $c\neq a$ and $b\neq d$.\par
\endproclaim
\vskip3mm A vertex $g$ in a directed graph $G$ is called a {\it
sink}, if the in-degree of $g$ is positive and the out-degree of
$g$ is zero. The dual concept of sink is called {\it source}.\par
\vskip3mm Let $Z_r(R)$ be the set of right zero divisors. In [6],
it is proved that for any ring $R$, $\Gamma(R)$ is connected if
and only if $Z_r(R)=Z_l(R)$, i.e., there exists no end-vertex
(sink or source) in $\Gamma(R)$. It is also proved that
$\Gamma(R)$ is connected for all artinian rings with two
sided-identity element. We now characterize all finite rings $R$
whose directed zero-divisor graph $\Gamma(R)$ is connected:\par
\proclaim{Theorem 2.4}
 For any finite ring $R$, the following statements
are equivalent:\par (1) The directed zero-divisor graph
$\Gamma(R)$ is connected;\par (2) Every one-sided identity element
of R is the two-sided identity of R; \par (2') Either R has
two-sided identity or R has no proper one-sided identity;\par (3)
There exists no end-vertex (sink or source) in
$\Gamma(R)$.([6])\par
\endproclaim
\vskip3mm
\demo{Proof}
$(1)\Longrightarrow(2)$ and $(3)\Longrightarrow(2).$
If $R$ has a left identity $e$ that is not a right
identity element of $R$, then there exists $a\in R$
such that $ae-a\not=0$.
In this case $ae-a\rightarrow b$ for all $b\in R^*$.
Then there exist at least two
left identity elements in $R$, say $e$ and $f=e+ae-a$.
Then $e$ and $f$ are sink vertices in $\Gamma(R)$. So
the directed graph $\Gamma(R)$ is not connected.\par

$(2)\Longrightarrow(3)$ and $(2)\Longrightarrow(1).$ Assume that
every one-sided identity element of R is the two-sided identity of
R. If $|R|<5$, then by Lemma 2.2, we know that there is no sink
(source) vertex in $\Gamma(R)$ and that $\Gamma(R)$ is connected.
In what follows, we assume $|R|\geq5$. Then by Lemma2.1, if
$a\rightarrow b$ for distinct vertices $a, b$, then there exist
$c,d\in R^*$ such that $c\not=a, d\not=b$ and in $\Gamma(R)$ there
is a walk $c\rightarrow a\rightarrow b\rightarrow d$. So there is
no sink or source vertex in $\Gamma(R)$. This proves (3). Now we
use the proof of corresponding result in [2] to finish our proof:
For any distinct $x,y\in Z(R)^*$, if $xy=0$, then $d(x,y)=1$. If
$xy\not=0$, then again by Lemma 2.1 and [6, Theorem 3.2], there
exists $a\neq x, b\neq y$ such that $xa=0=by$. If $a=b,$, then
$x\rightarrow a\rightarrow y$ and $d(x,y)=2$; If $a\neq b$ and
$ab=0$, then we have $x\rightarrow a\rightarrow b\rightarrow y$
and $d(x,y)\leq 3$; If $a\neq b$ and $ab\not=0$, then $ab\neq
x,ab\neq y$ since $xy\neq 0$, and there is a path $x\rightarrow ab
\rightarrow y$. In all cases, $d(x,y)\leq 3$.\hfill QED \par
\enddemo
\vskip3mm We remark that Theorem 2.4 actually holds for a still
wider class of rings - rings which are both right artinian and
left artinian, since Lemma 2.1 holds for artinian rings.
\par \vskip3mm
 Now let $R$ be a ring with proper left
identity element $e$. Denote
$$I_e=\{a\in R|ae=0\},R_e=\{a\in R|a=ae\}.$$ Then (1)
$I_e=ann_l(e)$ and $I_e$ is a two-sided ideal of $R$ with at leat
two elements;\par (2) $R_e$ is a subring of $R$ with identity
$e$;\par (3) $R=R_e\oplus I_e$ as left $R$-module;\par (4) There
is a ring isomorphism $R_e\cong R/I_e$.\par
 For any subset $M,N$ of $Z(R)^*$, let
$\Gamma(M, N)$ be the induced bipartite subgraph
of $\Gamma(R)$ and denote its edge set
as $E(M,N)$. For the graph $\Gamma(R)$, denote
its edge set as $E(R)$. If there is a cycle, we add 2
to the number of the edges of $\Gamma(R)$. Then the graph $\Gamma(R)$ is
completely determined by
the following four induced subgraphs: $\Gamma(R_e)$,
complete directed graph $\Gamma(I_e)$,
bipartite graphs $\Gamma(I_e^*, R_e^*)$ and
$\Gamma(R_e^*, R_e*\oplus I_e^*)$. We have\par
\proclaim{Proposition 2.5}
For a ring $R$ with proper left
identity element $e$, the graph $\Gamma(R)$ has
the following properties:\par
(1) For any $a\in I_e$, the out-degree of a
is $|R|+1$;\par
(2) The number of vertices of $\Gamma(R)$ is
$|R|-1=|R_e||I_e|-1$;\par
(3) $|E(R)|=|I|[|I|-1+E(R_e^*,I_e^* )+
E(R_e^*,R_e^*\oplus I_e^*)+(2-|I|)E(R_e)].$
\endproclaim
\vskip3mm
\demo{Proof}
(1) The number of directed edges from $I_e^*$ to
$R^*$ is $|R|(|I|-1)$;\par
(2) The number of directed edges from $R^*-I_e^*$ to
$I_e^*$ is
$$|I|[|E(R_e^*, I_e^*)|-(|I|-1)(|R_e|-1)]; $$
\indent (3) An element of $R^*-I_e^*$ has the form of
$a_i+x$, where $x=0$ or $x=b_j$. Besides,
$(a_i+x)(a_j+y)=a_i(a_j+y)$. Thus the number of
directed edges from $R^*-I_e^*$ to
$R^*-I_e^*$ is
$$|I|[|E(R_e)|+|E(R_e^*,R_e^*\oplus I_e^*)|-
(|I|-1)|E(R_e)|].$$
Note that $|R|=|R_e||I_e|$, then we
obtain our formula.\hfill QED
\enddemo
\vskip3mm
\proclaim{Proposition 2.6}
Let R be a finite ring with proper one-sided identity.
For distinct elements $x,y$ in $R$, either
$d(x,y)=\infty$ or $d(x,y)\leq 6$.\par
\endproclaim
\vskip3mm
\demo{Proof}
Let $R$ be a ring with proper left identity and
assume that
$e$ is a left identity of $R$. Then
$R^*=K\cup I_e^*\cup R_e^*$ and this is a disjoint union
where we assume
$$R_e^*=\{a_i|i=1,2,\cdots ,m\},
I_e^*=\{b_j|j=1,2,\cdots ,n\}$$ and
$K=I_e^*\oplus R_e^*$. For distinct $x,y\in R^*$,
assume that there exists a directed path from $x$ to $y$
in $\Gamma(R)$, say,
$x\to x_1\to x_2\to \cdots \to x_r\to y$. \par
(1) If $x\in I_e$, then $x\to z$ for any $z$. In this case,
$d(x,y)=1$.\par
(2) Assume $x\in R_e$. If $0\in xI_e^*$, then
we have a path $x\to I_e^*\to y$, $d(x,y)\leq 2$;
If $0\notin xI_e^*$, then $0\notin xK$. Without loss,
we can assume $x_i\in R_e$ for $i=1,2,3,4$. Since $R_e$ is a ring with
two-sided identity, $diam(R)\leq 3$. So there exists
a path $x\to x_1\to x_2\to x_3\to I_e^*\to y$. In this case,
$d(x,y)\leq 5$.\par
(3) The final case is $x\in K$.
Assume $x=a_1+b_1$. If $0\in xI_e^*$, then
$d(x,y)\leq 2$; If $0\notin xI_e^*$, then
$x_1\in R_e^*\cup K$. If $x_1\in K$, then we have
$0=(a_1+b_1)(a_s+b_t)=a_1a_s+a_1b_t$.
Then we have a path $x\to b_t\to y$, and $d(x,y)\leq 2$;
Hence we can assume $x_i\in R_e^*$ for $i=1,2,3,4$.
and $x_4\to I_e^*$. We have a path
$x\to x_1\to x_2\to x_3\to x_4 \to I_e^*\to y$ and hence
$d(x,y)\leq 6$. This completes the proof.\hfill QED
\enddemo
\vskip3mm
\proclaim{Corollary 2.7}
For any ring R with proper one-sided identity element
e, $R_e$ is a subring of $R$ and e is a two-sided
identity of $R_e$, and $diam(R)\leq 3+diam(R_e)$.
\endproclaim
\vskip3mm
We end this section with the following examples:\par
\noindent {\bf Example 2.8 }
For any field $F$, let $R$ be the $n$ by $n$ full
matrix ring over $F$ ($n>1$). Then the diameter of $\Gamma(R)$
is 2.\par
\indent We need only to prove the following facts:
for any $A,B\in Z(R)^*$, there exists $C\in Z(R)^*$
such that $AC=0$ and $CB=0$. First, there exist
invertible matrices $P,Q\in R$ such that the last
column of $AP$ is zero and the first row of $QB$ is
zero. Second, let
$$C=P
\pmatrix 0 &0 &\cdots &0\cr
\vdots &\vdots &\ddots &\vdots\cr
0 &0 &\cdots &0\cr
1 &0 &\cdots &0 \cr\endpmatrix
Q,$$
then we have $C\neq 0$,
$AC=0,CB=0$.\par
\vskip3mm
\noindent {\bf Example 2.9 }
Let $R$ be the $n$ by $n$ full matrix ring over
$\Bbb Z/2\Bbb Z$. Let $S$ be the non-unitary
subring of $R$ consisting of those matrices all
of whose rows are zero except the first row.
Then $Sink(S)=\{
\pmatrix 1 &\alpha \cr
0 &0 \cr
\endpmatrix\,|\,\alpha
 \},$
and $Sink(S)$ consists of all left identity elements
of $S$.
$I_e=
\{
\pmatrix 0 &\alpha \cr
0 &0 \cr
\endpmatrix\,|\,\alpha
 \}$
for $
e=
\pmatrix 1 &0 \cr
0 &0 \cr
\endpmatrix. $
$|S|=2^n,|I_e|=|Sink(S)|=2^{n-1}.$
$\Gamma(S)$ is a star-like directed graph
whose kernel is the complete graph
$K_{2^{n-1}-1}$, each vertex of which also connects to
$2^{n-1}$ sinks.\par
\vskip3mm
\noindent {\bf Example 2.10 }
Let $R=\Bbb Z/n\Bbb Z$ and denote
$S=
\{
\pmatrix a &b \cr
0 &0 \cr
\endpmatrix\,|\,a,b\in R
 \}$ .
Then $Sink(S)=
\{
\pmatrix a &b \cr
0 &0 \cr
\endpmatrix\,|\,a,b\in R, a=\overline{m}, (m,n)=1
 \}$.
$I_e=
\{
\pmatrix 0 &b \cr
0 &0 \cr
\endpmatrix\,|\,b\in R
 \}$
for
$e=
\pmatrix 1 &0 \cr
0 &0 \cr
\endpmatrix. $
Thus $|Sink(S)|=n\varphi (n)$, $|I_e|=n$, where
$\varphi (n)$ is the Eulerian number of $n$.
When $n\ge 3$, there is no source in $\Gamma(R)$.
The clique number of $\Gamma(R)$ is $n-1$.\par
\vskip3mm
\head 3. Sinks and Sources of finite rings\endhead
\vskip3mm
\proclaim{Proposition 3.1}
(1) For any ring R, if there exists a source b
in $\Gamma(R)$ with $b^2=0$, then
$R=\{0,a,b,c\}$, where $a$ and $c$ are left
identity elements, $ba=0=bc$;  \par
(2) For any ring R, if there exists a sink b
in $\Gamma(R)$ with $b^2=0$, then
$R=\{0,a,b,c\}$, where $a$ and $c$ are
right identity elements, $ab=0=cb$.
\endproclaim
\vskip3mm
\demo{Proof}
(1) Let $b$ be a source
in $\Gamma(R)$ with $b^2=0$. Then consider
the left $R$-module epimorphism
$\eta : R\to Rb, r\mapsto rb$. Notice that $ker(\eta )=ann_l(b)=\{0,b\} $
is a
submodule of $R$. For any $r\in R-\{0,b\}$, we have
$rb=b$. We conclude $|R|\leq 4$, since otherwise,
let $a,c,d\in R-\{0,b\}$ be
distinct elements. Then we have $(a-c)b=b-b=0=(a-d)b$.
Therefore, $a-c=a-d$ and thus $c=d$, a contradiction.
\par
Now by the proof of Lemma 2.2 and Theorem 2.4, $R$ has proper left
identity element. So $R$ contains at least two left identity
element, say $a$ and $c$. Thus $R=\{0,a,b,c\}$, where $a$ and $c$
are left identity elements, $ba=0=bc$. In this case, we have
$a\leftarrow b\leftrightarrows b \rightarrow c$. \par (2) The
proof is dual to the above proof.
\enddemo
\vskip3mm

The structure of rings $R$ with $|R|\leq 4$, as well as the
related graph $\Gamma(R)$, is rather clear (see Lemma 2.2 and the
listing following Lemma 2.2). As for finite rings $R$ with
$|R|\geq 5$, we have the following:\par \vskip3mm
\proclaim{Proposition 3.2} Let $R$ be a finite ring with proper
left identity elements. If $|R|\geq 5$, then\par (1) $\Gamma(R)$
contains at least two sinks;\par (2) For any sink $r$ in
$\Gamma(R)$, $r^2\neq 0$;\par (3) $\Gamma(R)$ contains no sources.
\endproclaim
\vskip3mm
\demo{Proof}
Since every left identity element of $R$ is a sink in
$\Gamma(R)$, thus $\Gamma(R)$ has at least two sinks.
Let $e$ be any left identity element of $R$.
Then $I_e^*$ is not empty, and for any
$x\in I_e^*,y\in R^*$,
we have $xy=0$. Thus every non-zero element of $R$ is
a vertex of $\Gamma(R)$ and, if $\Gamma(R)$ has source,
the source vertex must lies in $I_e$. Of course,
$\Gamma(R)$ has source
vertex if and only if $|I_e|=2$ and there is no directed
edge from $R^*-I_e^*$ to $I_e^*$. This is equivalent to
saying that for the left identity element $e$, there exists
a non-zero element
$b$ of $R$, such that $ann_l(e)=\{0,b\}$ and
$ann_l(b)=\{0,b\}$. So
if $\Gamma(R)$ has a
source, this source is $b$. Then we have $b^2=0$,
contradicting with the assumption and Proposition 3.1.
With the same reason, sinks $r\in \Gamma(R)$ must satisfy
$r^2\neq 0$. \hfill QED
\enddemo
\vskip3mm
\proclaim{Corollary 3.3}
Let $R$ be a finite ring with proper left identity
elements. Assume that  $R$ has at least five elements.
For any left identity $e$ of $R$, if
$ann_l(e)=\{0,b\}$ for some non-zero
element $b$, then in $\Gamma(R)$, the out-degree of
$b $ is $|R|-1$ and the in-degree of $b$ is positive.
\endproclaim
\vskip3mm
Similarly, we have\par
\proclaim{Proposition 3.4}
Let $R$ be a finite ring with proper right identity
elements. If $R$ has at least five elements, then\par
(1) $\Gamma(R)$ contains at least two sources;\par
(2) For any source $r$ in $\Gamma(R)$, $r^2\neq 0$;\par
(3) $\Gamma(R)$ contains no sink.
\endproclaim
\vskip3mm
\proclaim{Corollary 3.5}
Let $R$ be a finite ring with proper right identity
elements. If $R$ has at least five elements,
then for any right identity $e$ of $R$ with
$ann_l(e)=\{0,b\}$ for some non-zero
element $b$, the in-degree of
$b$ in $\Gamma(R)$ is $|R|-1$ and the out-degree of $b$ is positive.
\endproclaim
\vskip3mm
Recall that a {\it network} $N$ is a directed graph with
exactly
one sink vertex $k$ and a unique source vertex $c$ such that
$c$ connects to every vertex of $N$ and
every vertex of $N$ connects to $k$. By Lemma 2.2,
Theorem 2.5,
Propositions 3.1,3.5 and 3.7, we immediately have,\par
\proclaim{Corollary 3.6} For any finite ring R
with at least five elements,\par
(1) $\Gamma(R)$ can not contain sink and source at
the same time;\par
(2)([7, Cor. 3.2])
 $\Gamma(R)$ is not a network for any finite ring.
\endproclaim
\vskip3mm
\demo{Proof}
(2) is an obvious consequence of (1).
To prove (1), we list all our previous results on rings
with proper one-sided identity in a
single place as follows:\par
(a) If $|R|\leq 4$, then $\Gamma(R)$ is one of the following
$$\circ \leftarrow \circ \rightarrow \circ,
\circ\rightarrow\circ \leftarrow\circ.$$
\indent (b) If $|R|\geq 5$ and $R$ has proper left identity elements,
then $\Gamma(R)$ contains at least two sinks,
but in it there is no source. \par
(c) If $|R|\geq 5$ and $R$ has proper right identity elements,
then $\Gamma(R)$ contains at least two sources,
but in it there is no sink.
\enddemo
\vskip3mm
\noindent {\bf Definition 3.7}
Suppose that $R$ has proper left (respectively, right)
identity element.
An element $r\in R$ is called {\bf strongly right
(left) invertible}, if for {\it any} left (right) identity
element $e$ of $R$, $r$ has unique right (left) inverse
$s$ relative to this $e$. In this case, such $s$ has
more than one
left (right) inverse relative to the same left
identity $e$. Obviously, $r$ is strongly right invertible
if for {\it some} left identity
element $e$ of $R$, $r$ has a unique right inverse
$s$ relative to this $e$.\par
\vskip3mm
\proclaim{Proposition 3.8}
For an element $r$ in a finite ring $R$ satisfying
$r^2\neq 0$, $r$ is a sink vertex (source vertex) in $\Gamma(R)$,
if and only if $r$ is strongly right (respectively, left)
invertible in $R$.
\endproclaim
\vskip3mm
\demo{Proof}
(1) Suppose $r$ is a sink vertex of $\Gamma(R)$. We have
$ry\neq 0$ for all $y\in R^*$ since $r^2\neq 0$.
Since $R$ is a finite ring, we have $rR=R$. Assume
$re=r$. Then from $res=rs\neq 0$ for all $s\in R^*$,
we obtain $es=s$ for all $s\in R$. Thus $e$ is a left
identity of $R$. For any left identity $f\in R$,
let $ru=f$. Then $u$ is unique relative to the $f$,
since $r$ is a sink vertex. For the same reason,
there is $r\neq v\in R^*$
such that $vr=0$. Hence $0=vru=vf$. Thus $f$ is not a
right identity element of $R$. In this case, we have
$(r+v)u=ru+v(fu)=f$, where $r+v\neq r$. Thus this
$u$ has
at least two left inverses (say, $r$ and $r+v$) relative
to the left identity $f$. Thus $r$ is strongly right
invertible in $R$. Notice that $R$ has at least two
sink vertex in this case.\par
(2) Conversely, Suppose that $r$ is strongly right
invertible in $R$. Then by definition, there is left,
but not right,
identity element in $R$, and for any left identity
$e\in R$, $r$ has a unique right inverse
$u$ relative to this $e$. If there is a path
$r\rightarrow v$, then we have $r(u+v)=e$ and this implies
$v+u=u$, a contradiction. Finally, suppose $x=ae-a\neq 0$,
then $x\neq r$ and $xr=0$. Thus $r$ is a sink vertex of
$\Gamma(R)$.\par
The proof of the other case is similar.\hfill QED \par
\enddemo

\proclaim{Corollary 3.9}Let R be a finite ring.\par
 (1)If $\Gamma(R)$ contains exactly one source
(respectively, sink),
then $|R|=4$, and $R^*$ has the form of
$$a\leftarrow b\rightarrow b\rightarrow c;$$
$$(respectively, a\rightarrow b\rightarrow b
\leftarrow c);$$ \par
(2)Let $|R|\geq 5$. Then an element $r$ of $R^*$ is a
sink (source) if and  only if $r$ is strongly right (left)
invertible in $R$. In this case, $r^2\neq 0$ and
in $\Gamma(R)$ there are at least two sinks (sources).\par
\endproclaim
\vskip3mm
\head 4. Sink(R), Sour(R) and network
\endhead
\vskip3mm
All rings in this section have at least five
elements.\par
\noindent {\bf Definition 4.1} For any ring $R$,
denote \par
$Sink(R)$ = $\{$sinks in $\Gamma(R)$ $ \}$,
$Sour(R)$ = $\{$sources in $\Gamma(R)$ $ \},$ \par
$Inv_r(R)$ = $\{$ strongly right  invertible
elements of
$R$ relative to some proper left identity $ \}$,  \par
$Inv_l(R)$ = $\{$  strongly left invertible
elements of $R$ relative to some proper right identity$ \}$.
\vskip3mm
\proclaim{Proposition 4.2}
Let R be any ring with at least five elements.\par
(1) If $Sink(R)$ ( respectively, $Sour(R)$ ) is not empty,
then it is a left ( right ) cancellative multiplicative semigroup;\par
(2) If $Inv_r(R)$ (  $Inv_l(R)$  ) is not empty, then
it is also a
left ( right ) cancellative multiplicative semigroup and
$Inv_r(R)\subseteq Sink(R)$
( $Inv_l(R)\subseteq Sour(R)$ );\par
(3) $Sink(R)=Z_r(R)-Z_l(R)$, $Sour(R)=Z_l(R)-Z_r(R);$\par
(4) $Z(R)^*$ has a disjoint decomposition
$$Z(R)^*=Sour(R)\cup (Z_r(R)\cap Z_l(R))\cup Sink(R).$$
\endproclaim
\demo{Proof}
(1) Assume that $Sink(R)$ is nonempty. For any
$a,b\in Sink(R)$,
there exists $x\in R^*$ such that $x\neq a$ and
$xa=0$. By Proposition 3.1, $ann_r(a)=0$.
So $Sink(R)$ is left cancellative and
$Sink(R)=Z_r(R)-Z_l(R)$. We assert
$x\neq ab$, since otherwise,
we would have $abx=0$, which implies $x=0$. Thus
$ab\in Sink(R)$, since $xab=0$. Hence $Sink(R)$ is
a semigroup under the multiplication of $R$.\par
(2) Suppose that $Inv_r(R)$ is not empty.
For any proper left identity
elements $e,f$ of $R$ and
any $a,b\in Inv_r(R)$, let $ax=e,f=by$. Then
$(ab)(yx)=afx=ax=e$. Thus $ab$ is right invertible.
If $(ab)c=e=(ab)d$, then $bc=bd$, and $c=d$ since
$Inv_r(R)\subseteq Sink(R)$. Thus $ab$ is strongly
right invertible, i.e., $ab\in Inv_r(R)$ and
hence $Inv_r(R)$ is a
semigroup. \par
The other case is dual to the sink case.
\hfill QED
\enddemo
\vskip3mm
{\bf Remark } For some fixed proper left identity
$e$, let
$Inv_{r_e}^{-1}(R)$ =
$\{$$u\in R|au=e$ for some $a\in Inv_r(R)$ $ \}$.
Then it is easy to verify that $Inv_{r_e}^{-1}(R)$ is a
multiplicative semigroup with identity $e$. \par
\vskip3mm
\proclaim{Proposition 4.3}
For any ring $R$ with at least five
elements, R contains a proper left identity
if and only if the following two conditions hold:
(1) $Sink(R)\neq \varnothing$; and (2) There exists an
$x\in Sink(R)$ such that $x[Sink(R)]=Sink(R)$.\par
In this case, $\Gamma(R)$ contains no source while
$Sink(R)$ contains at least two elements.
\endproclaim
\vskip3mm
\demo{Proof}
If $R$ contains a proper left identity $e$,
then $e\in Sink(R)$ and $e[Sink(R)]=Sink(R)$.\par
Conversely, assume that $Sink(R)\neq \varnothing$ and
$x[Sink(R)]=Sink(R)$ for some $x\in Sink(R)$.
Let $xe=x$, $e\in Sink(R)$. Then $x(ey-y)=0$ for all
$y\in R^*$. Then by Proposition 3.1, $ey=y$ and hence,
$e$ is a left identity of $R$. Since $e$
is a zero-divisor, it is a proper left identity.\par
If there is a source in $\Gamma(R)$, then it must
lie in $I_e$. Then $I_e=\{0,b\}$ for some nonzero
$b$. Then $b$ is a source of $\Gamma(R)$ with
$b^2=0$. Then by Proposition 3.1, $|R|=4$ and there
is no sink in $\Gamma(R)$, a contradiction.
So in this case, there is no source in $\Gamma(R)$.
\hfill QED
\enddemo
\vskip3mm
Since $Sink(R)$ is left cancellative, we immediately
have\par
\proclaim{Corollary 4.4}
For any ring R with at least five
elements, if
$0<|Sink(R)|<\infty$, then $Sink(R)=Inv_r(R)$ and
there is no source in $\Gamma(R)$.
\endproclaim
\vskip3mm
The following two results are duals of 4.3 and 4.4\par
\proclaim{Proposition 4.5}
For any ring $R$ with at least five
elements, R contains a proper right identity,
if and only if the following two conditions hold:
(1) $Sour(R)\neq \varnothing$; and (2) There exists an
$y\in Sour(R)$ such that $[Sour(R)]y=Sour(R)$.\par
In this case, $\Gamma(R)$ contains no sink
but it contains at least two sources.
\endproclaim
\vskip3mm
\proclaim{Corollary 4.6}
For any ring R with at least five
elements, if $0<|Sour(R)|<\infty$, then $Sour(R)=Inv_l(R)$ and
there is no sink in $\Gamma(R)$.
\endproclaim
\vskip3mm
As a combination of Propositions 4.2, 4.3 and 4.5, we have\par
\proclaim{Proposition 4.7}
For any ring R, one-sided identity of R are two-sided
identity if and only if $\Gamma(R)$ satisfies one of the
following conditions:\par
(1) $Sink(R)=\varnothing$ and $Sour(R)=\varnothing$;\par
(2) $Sink(R)=\varnothing$, $Sour(R)\neq\varnothing$
and for any $t\in Sour(R)$,
$[Sour(R)]t\subset Sour(R)$; (In this case,
$\Gamma(R)$ has infinitely many sources);\par
(3) $Sour(R)=\varnothing$, $Sink(R)\neq\varnothing$
and for any $s\in Sink(R)$,
$s[Sink(R)]\subset Sink(R)$; (In this case,
$\Gamma(R)$ has infinitely many sinks);\par
(4) $Sink(R)\neq\varnothing$
and for any $s\in Sink(R)$, $s[Sink(R)]\subset Sink(R)$.
At the same time, $Sour(R)\neq\varnothing$
and for any $t\in Sour(R)$,
$[Sour(R)]t\subset Sour(R)$. ( In this case, $\Gamma(R)$
has infinitely many sinks and infinitely many sources.)
\endproclaim
\vskip3mm \proclaim{Corollary 4.8} Suppose that in a ring R,
one-sided identity element is two-sided identity. If in addition,
R satisfies descending chain condition on principal left
(respectively, right) ideals, then $\Gamma(R)$ contains no source
(sink). In particular, if R is a left and right artinian ring, and
one-sided identity element in R is two-sided identity, then
$\Gamma(R)$ contains neither source no sink.\par
\endproclaim
\vskip3mm
\demo{Proof}
It is easy to verify that for $K=Sink(R)$ and any
$t\in K$, $tK\subset K$ if and only if for any
(or some) $s\in K$, $stK\subset sK$, if and only if
$stR\subset sR$. So, if $R$ has DCC on right
principal ideals, then in Proposition 4.7, cases (3)
and (4) could not occur. Thus $\Gamma(R)$ contains no
sinks.\hfill QED
\enddemo
\vskip3mm
Finally, as a corollary of Propositions 3.1, 4.3,
4.5 and 4.7, we have\par
\proclaim{Corollary 4.9}
For any ring R, $\Gamma(R)$ is not a network.
\endproclaim
\vskip5mm \centerline{\bf REFERENCES}
 \vskip3mm
\ref\no 1
\by D. D. Anderson and M. Naseer
\paper Beck's coloring of a commutative ring
\jour J. Algebra
\vol 159
\yr 1993
\pages  500-514
\endref

\ref\no 2
\by D. F. Anderson and P. S. Livingston
\paper The zero-divisor graph of a commutative ring
\jour  J. Algebra
\vol 217
\yr 1999
\pages 434-447
\endref

\ref\no 3
\by  D. F. Anderson, Ron Levy, J. Shapiro
\paper Zero-divisor graphs, von Neumann regular rings, and Boolean algebras
\jour J. Pure Applied Algebra
\vol 180
\yr 2003
\pages 221-241
\endref

\ref\no 4
\by I. Beck
\paper Coloring of commutative rings
\jour J. Algebra
\vol 116
\yr 1988
\pages 208-226
\endref

\ref\no 5
\by F. R. DeMeyer, T. McKenzie, and K. Schneider
\paper The zero-divisor graph of a commutative semigroup
\jour Semigroup Forum
\vol 65
\yr  2002
\pages  206-214
\endref

\ref\no 6
\by S. P. Redmond
\paper The zero-divisor graph of a non-commutative ring
\jour Internat. J. Commutative  Rings
\vol 1(4)
\yr 2002
\pages 203-211
\endref

\ref\no 7
\by S. P. Redmond
\paper Structure in the zero-divisor graph of a
non-commutative ring
\jour Preprint
\vol
\yr 2003
\pages
\endref

\ref\no 8
\by R. J. Wilson
\book  Introduction to Graph Theory
\publ  Longman Inc.
\publaddr New York
\yr Third Edition 1985
\endref

\enddocument